\documentclass{amsart}
\usepackage{latexsym}
\usepackage{amsmath,amsthm,amsfonts,amscd,eucal}
\usepackage{graphicx}
\usepackage{epsfig}

\numberwithin{equation}{section}

\def\ca{{\mathcal A}}
\def\cb{{\mathcal B}}
\def\cc{{\mathcal C}}

\def\cf{{\mathcal F}}

\def\cai{{\mathcal I}}

\def\cl{{\mathcal L}}
\def\cam{{\mathcal M}}
\def\cn{{\mathcal N}}

\def\cp{{\mathcal P}}

\def\car{{\mathcal R}}
\def\cs{{\mathcal S}}

\def\bc{{\mathbb C}}

\def\bn{{\mathbb N}}
\def\br{{\mathbb R}}
\def\bz{{\mathbb Z}}

\def\a{\alpha}

\def\g{\gamma}        \def\G{\Gamma}
\def\d{\delta}        \def\D{\Delta}
\def\eps{\varepsilon}

\def\th{\vartheta}

\def\l{\lambda}       \def\La{\Lambda}
\def\m{\mu}

\def\r{\rho}
\def\s{\sigma}        
    
\def\t{\tau}

\def\o{\omega}        \def\O{\Omega}

\def\itm#1{\item{$(#1)$}}
\def\ov{\overline}
\newcommand{\vol}{\text{vol}}
\newcommand{\set}[1]{\left\{#1\right\}}
\def\de{\partial}
\newcommand{\conv}{\text{conv}\,}

\newcommand{\Det}{\text{det}}
\newcommand{\ee}{\tilde{e}}

\newcommand{\Ci}{\cc}  
\renewcommand{\Re}{\car}  
\renewcommand{\Pr}{\cp}  

\def\square{\hbox{$\sqcap\!\!\!\!\sqcup$}}
\def\endproof{\hfill\square}

\newtheorem{Thm}{Theorem}[section]
\newtheorem{Cor}[Thm]{Corollary}
\newtheorem{Prop}[Thm]{Proposition}
\newtheorem{Lemma}[Thm]{Lemma}

\theoremstyle{definition}
\newtheorem{Dfn}[Thm]{Definition}
\newtheorem{exmp}[Thm]{Example}

\theoremstyle{remark}
\newtheorem{rem}[Thm]{Remark} 
\newtheorem{ack}{Acknowledgement} 

\begin{document}

\title[Ihara's zeta function for periodic graphs and
its approximation]{Ihara's zeta function for periodic graphs\\ and
its approximation in the amenable case}
\author{Daniele Guido, Tommaso Isola, Michel L. Lapidus}%
\date{August 9, 2006}%
\address{(D.G., T.I.) Dipartimento di Matematica, Universit\`a di Roma ``Tor
Vergata'', I--00133 Roma, Italy.}%
\email{guido@mat.uniroma2.it, isola@mat.uniroma2.it}%
\address{(M.L.L.) Department of Mathematics, University of California,
Riverside, CA 92521-0135, USA.}%
\email{lapidus@math.ucr.edu}%
\thanks{The first and second authors were partially
supported by MIUR, GNAMPA and by
the European Network ``Quantum Spaces - Noncommutative Geometry"
HPRN-CT-2002-00280. The third author was partially supported
by the National Science Foundation, the Academic Senate of the 
University of California, and GNAMPA}%
\subjclass[2000]{Primary 05C25,11M41, 46Lxx;
Secondary 05C38, 11M36, 30D05. }%
\keywords{Periodic graphs, Ihara zeta function,
    analytic determinant, determinant formula, functional equations,
    amenable groups, amenable graphs, approximation by finite
    graphs.}%

\maketitle
\bigskip

\begin{abstract}
    In this paper, we give a more direct proof of the results by Clair
    and Mokhtari-Sharghi \cite{ClMS1} on the zeta functions of
    periodic graphs.  In particular, using appropriate
    operator-algebraic techniques, we establish a determinant formula
    in this context and examine its consequences for the Ihara zeta
    function.  Moreover, we answer in the affirmative one of the
    questions raised in \cite{GrZu} by Grigorchuk and
    $\dot{\text{Z}}$uk.  Accordingly, we show that the zeta function
    of a periodic graph with an amenable group action is the limit of
    the zeta functions of a suitable sequence of finite subgraphs.    
\end{abstract}



\setcounter{section}{-1}

\section{Introduction}

 The zeta functions associated to finite graphs by Ihara \cite{Ihara},
 Hashimoto \cite{HaHo,Hashi}, Bass \cite{Bass} and others, combine
 features of Riemann's zeta function, Artin L-functions, and Selberg's
 zeta function, and may be viewed as analogues of the Dedekind zeta
 functions of a number field.  They are defined by an Euler product
 and have an analytic continuation to a meromorphic function
 satisfying a functional equation.  They can be expressed as the determinant
 of a perturbation of the graph Laplacian and, for  Ramanujan graphs,
 satisfy a counterpart of the Riemann hypothesis \cite{StTe}.  Other relevant
 papers are \cite{Sunada,Hashi1,Hashi2,Serre1,North,FoZe,KoSu,StTe2,StTe3,HSTe,Bartholdi,MiSa}.

 In differential geometry, researchers have first studied compact
 manifolds, then infinite covers of those, and finally, noncompact
 manifolds with greater complexity.  Likewise, in the graph setting,
 one passes from finite graphs to infinite periodic graphs, and then
 possibly to other types of infinite graphs.  In fact, the definition
 of the Ihara zeta function was extended to (countable) periodic
 graphs  by Clair and Mokhtari-Sharghi
 \cite{ClMS1}, and a corresponding determinant formula was proved.
 They deduce this result as a specialization of the treatment of
 group actions on trees (the so-called theory of tree lattices, as
 developed by Bass, Lubotzky and others, see \cite{BaLu}).  We mention
 \cite{GILa02} for a recent review of some results on zeta functions
 for finite or periodic simple graphs, and \cite{GrZu,ClMS1,ClMS2,Clair}
 for the computation of the Ihara zeta function of several periodic simple graphs.

 In \cite{GrZu}, Grigorchuk and $\dot{\text{Z}}$uk defined zeta
 functions of infinite discrete groups, and of some class of infinite
 periodic graphs (which they call residually finite), and asked how to
 obtain the zeta function of a periodic graph by means of the zeta
 functions of approximating finite subgraphs, in the case of amenable
 or residually finite group actions.

 The purpose of the present work is twofold: first, to give a different proof of the main result
 obtained by Clair and Mokhtari-Sharghi in \cite{ClMS1}; second,  to answer in the affirmative
 one of the questions raised by Grigorchuk and $\dot{\text{Z}}$uk in \cite{GrZu}.

 As for the first point, some combinatorial results in Section \ref{sec:Prelim}
 give a more direct proof of the determinant formula in Theorem \ref{Thm:det.formula}. Moreover,
 the theory of analytic determinants developed in Section \ref{sec:AnalyticDet} allows us to use
 analytic functions instead of formal power series in that formula, as well as to establish functional
 equations for suitable completions of the Ihara zeta function, generalizing results contained in \cite{GILa02}.

 As for the second point, we take advantage of the technical framework developed in this paper to
 show,  in the case of amenable group actions,
 that the Ihara zeta function is indeed the limit of the zeta
 functions of a suitable sequence of approximating finite graphs.  For
 the sake of completeness, we mention that, in \cite{ClMS2}, Clair and
 Mokhtari-Sharghi have given a positive answer in the case of
 residually finite group actions.

 This paper is organized as follows. We start in Section \ref{sec:Prelim} by recalling some notions from
 graph theory and prove all the combinatorial results we need in the
 following sections.  In Section \ref{sec:Zeta}, we then define the
 analogue of the Ihara zeta function and show that it is a
 holomorphic function in a suitable disc, while, in Section
 \ref{sec:DetFormula}, we prove a corresponding determinant formula,
 which relates the zeta function with the Laplacian of the graph.  The
 formulation and proof of this formula requires some care because it
 involves the definition and properties of a determinant for bounded
 operators (acting on an infinite dimensional Hilbert space and)
 belonging to a von Neumann algebra with a finite trace.  This
 issue is addressed in Section \ref{sec:AnalyticDet}.  In Section
 \ref{sec:functEq}, we establish several functional equations for
 various possible completions of the zeta function.  In the final
 section, we prove the approximation result mentioned above.

 In closing this introduction, we note that  in \cite{GILa01}
 we define and study the Ihara zeta functions attached to a
 new class of infinite graphs, called self-similar fractal graphs,
 which have greater complexity than the periodic ones.

 The contents of this paper have been presented  at the $21^{st}$ conference on Operator Theory
in Timisoara (Romania) in July 2006.

\section{Preliminary results}\label{sec:Prelim}

 We recall some notions from graph theory, following \cite{Serre}.  A
 {\it graph} $X=(VX,EX)$ consists of a collection $VX$ of objects,
 called {\it vertices}, and a collection $EX$ of objects called
 (oriented) {\it edges}, together with two maps $e\in EX\mapsto
 (o(e),t(e))\in VX\times VX$ and $e\in EX\mapsto \ov{e}\in EX$,
 satisfying the following conditions: $ \ov{\ov{e}}=e$,
 $o(\ov{e})=t(e)$, $\forall e\in EX$.  The vertex $o(e)$ is called the
 {\it origin} of $e$, while $t(e)$ is called the {\it terminus} of
 $e$.  The edge $e$ is said to join the vertices $u:=o(e)$, $v:=t(e)$,
 while $u$ and $v$ are said to be {\it adjacent}, which is denoted
 $u\sim v$.  The edge $e$ is called a {\it loop} if $o(e)=t(e)$.  The
 {\it degree} of a vertex $v$ is $\deg(v):= |\set{e\in EX: o(e)=v}|$,
 where $|\cdot|$ denotes the cardinality.  A {\it path} of length $m$
 in $X$ from $u=o(e_{1})\in VX$ to $v=t(e_{m})\in VX$ is a sequence of
 $m$ edges $(e_{1},\ldots,e_{m})$, where $o(e_{i+1})= t(e_i)$, for
 $i=1,...,m-1$.  In the following, the length of a path $C$ is denoted
 by $|C|$.  A path is {\it closed} if $u=v$.  A graph is said to be
 {\it connected} if there is a path between any pair of distinct
 vertices.

 The couple $\set{e,\ov{e}}$ is called a {\it geometric edge}.  An
 {\it orientation} of $X$ is the choice of one oriented edge for each
 couple, which is called {\it positively oriented}.  Denote by $E^{+}X$ the
 set of  positively oriented edges.  Then the other edge of each couple
 will be called negatively oriented, and denoted $\ov{e}$, if $e\in
 E^{+}X$.  The set of negatively oriented edges is denoted $E^{-}X$.
 Then $EX = E^{+}X \cup E^{-}X$.

 In this paper, we assume that the graph $X=(VX,EX)$ is connected,
 countable [$i.e.$ $VX$ and $EX$ are countable sets] and with bounded
 degree [$i.e.$ $d:= \sup_{v\in VX}\deg(v) <\infty$].  We also choose,
 once and for all, an orientation of $X$.

 Let $\G$ be a countable discrete subgroup of automorphisms of $X$,
 which acts
 \begin{enumerate}
    \item without inversions, $i.e.$ $\g(e)\neq \ov{e}, \forall
    \g\in\G, e\in EX$,

    \item discretely, $i.e.$ $\G_{v}:=\set{\g\in\G: \g v=v}$ is
    finite, $\forall v\in VX$,

    \item with bounded covolume, $i.e.$
    $\displaystyle{
    \vol(X/\G) := \sum_{v\in \cf_{0}} \frac{1}{|\G_{v}|} <\infty,
    }$
    where $\cf_0 \subset VX$ contains exactly one representative for each equivalence class in $VX/\G$.
 \end{enumerate}

 \noindent We note that the above bounded covolume property is equivalent to
 $$
 \vol(EX/\G) := \sum_{e\in \cf_{1}} \frac{1}{|\G_{e}|} <\infty,
 $$
 where $\cf_1 \subset EX$  contains exactly one representative for each equivalence class in $EX/\G$.

 Let us now define two useful unitary representations of $\G$.

 Denote by $\ell^{2}(VX)$ the Hilbert space of functions $f:VX\to\bc$
 such that $\|f\|^{2} := \sum_{v\in VX} |f(v)|^{2} <\infty$.  A
 unitary representation of $\G$ on $\ell^{2}(VX)$ is given by
 $(\l_{0}(\g)f)(x):= f(\g^{-1}x)$, for $\g\in\G$, $f\in\ell^{2}(VX)$,
 $x\in VX$.  Then the von Neumann algebra $\cn_{0}(X,\G):= \{
 \l_{0}(\g) : \g\in\G\}'$ of all the bounded operators on
 $\ell^{2}(VX)$ commuting with the action of $\G$, inherits a trace
 given by
 \begin{equation}\label{eq:TraceOnVertices}
     Tr_{\G}(A) := \sum_{x\in \cf_0} \frac{1}{|\G_{x}|}A(x,x), \
     A\in\cn_{0}(X,\G).
 \end{equation}

 Analogously, denote by $\ell^{2}(EX)$ the Hilbert space of functions
 $\o:EX\to\bc$ such that $\| \o \|^{2} := \sum_{e\in EX} |\o(e)|^{2}
 <\infty$.  A unitary representation of $\G$ on $\ell^{2}(EX)$ is
 given by $(\l_{1}(\g)\o)(e):= \o(\g^{-1}e)$, for $\g\in\G$,
 $\o\in\ell^{2}(EX)$, $e\in EX$.  Then the von Neumann algebra
 $\cn_{1}(X,\G):= \{ \l_{1}(\g) : \g\in\G\}'$ of all the bounded
 operators on $\ell^{2}(EX)$ commuting with the action of $\G$,
 inherits a trace given by
 \begin{equation}\label{eq:TraceOnEdges}
     Tr_{\G}(A) := \sum_{e\in \cf_1} \frac{1}{|\G_{e}|}A(e,e), \
     A\in\cn_{1}(X,\G).
 \end{equation}

 At this stage, we need to introduce some additional terminology from
 graph theory.

 \begin{Dfn}[Reduced Paths]
    \itm{i} A path $(e_{1},\ldots,e_{m})$ has {\it backtracking} if
    $e_{i+1}=\ov{e}_{i}$, for some $i\in\{1,\ldots,m-1\}$.  A path
    with no backtracking is also called {\it proper}.

    \itm{ii} A  closed path is called {\it primitive} if it is not
    obtained by going $n\geq 2$ times around some other closed path.

    \itm{iii} A proper closed path $C=(e_{1},\ldots,e_{m})$ has a {\it
    tail} if there is $k\in\bn$ such that $e_{m-j+1} = \ov{e}_{j}$,
    for $j=1,\ldots,k$.  Denote by $\Ci$ the set of proper tail-less
    closed paths, also called {\it reduced} closed paths.
 \end{Dfn}

 \begin{Dfn}[Cycles]
    Given closed paths $C=(e_{1},\ldots,e_{m})$,
    $D=(e'_{1},\ldots,e'_{m})$, we say that $C$ and $D$ are {\it
    equivalent}, and write $C\sim_{o} D$, if there is $k\in\bn$ such that
    $e'_{j}=e_{j+k}$, for all $j$, where $e_{m+i}:=e_i$, that is, the
    origin of $D$ is shifted $k$ steps with respect to the origin of
    $C$.  The equivalence class of $C$ is denoted $[C]_o$.  An
    equivalence class is also called a {\it cycle}.  Therefore, a closed
    path is just a cycle with a specified origin.

    Denote by $\Re$ the set of reduced cycles, and by $\Pr\subset\Re$
    the subset of primitive reduced cycles, also called {\it prime cycles}.
 \end{Dfn}

 \begin{Dfn}[Equivalence relation]
    \itm{i} Given $C$, $D\in\Ci$, we say that $C$ and $D$ are
    $\G$-{\it equivalent}, and write $C \sim_{\G} D$, if there is an
    isomorphism $\g\in \G$ such that $D=\g(C)$.  We denote by
    $[\Ci]_{\G}$ the set of $\G$-equivalence classes of reduced closed
    paths.

    \itm{ii} Similarly, given $C$, $D\in\Re$, we say that $C$ and $D$
    are $\G$-{\it equivalent}, and write $C \sim_{\G} D$, if there is
    an isomorphism $\g\in \G$ such that $D=\g(C)$.  We denote by
    $[\Re]_{\G}$ the set of $\G$-equivalence classes of reduced
    cycles, and analogously for the subset $\Pr$.
 \end{Dfn}

 \begin{rem}
     In the rest of the paper, we denote by $\Ci_{m}$ the subset of
     $\Ci$ consisting of closed paths of length $m$.  An analogous
     meaning is attached to $\Re_{m}$ and $\Pr_{m}$.
 \end{rem}

 Our proof of formula $(iv)$ in Theorem \ref{Prop:power.series} requires
 a generalization of a result by Kotani and Sunada \cite{KoSu} to infinite covering graphs.
 This is done in Proposition \ref{Lem:T}, whose proof depends on a new combinatorial result contained in Lemma \ref{counting}.

 Define the {\it effective
 length} of a cycle $C$, denoted by $\ell(C)$, as the length of the prime cycle
 underlying $C$, and observe that $\ell(C)$ is constant on the
 $\G$-equivalence class of $C$.  Therefore, if $\xi\in[\Re]_{\G}$, we
 can define $\ell(\xi):= \ell(C)$, for any representative $C\in\xi$.  Recall that, for any cycle $C$, the {\it stabilizer} of $C$ in $\G$ is the
 subgroup $\G_{C}:= \{\g\in\G : \g(C)=C\}$.  
 Moreover, if $C_{1},\, C_{2}\in \xi$, then the stabilizers
 $\G_{C_{1}},\, \G_{C_{2}}$ are conjugate subgroups in $\G$, and we
 denote by $\cs(\xi)$ their common cardinality. 


 For the purposes of the next few results, for any closed path
 $D=(e_{0},\ldots,e_{m-1})$, we also denote $e_{j}$ by $e_{j}(D)$.


 \begin{Lemma}\label{counting}
     Let $\xi\in[\Re_{m}]_{\G}$. Then
     $$
     \sum_{e\in\cf_{1}} \frac{1}{|\G_{e}|} |\set{D\in\Ci_{m}:
     [D]_{o,\G}=\xi, e_{0}(D)=e}| = \frac{\ell(\xi)}{\cs(\xi)}.
     $$
 \end{Lemma}
 \begin{proof}
     Let us first observe that, if $C_{1},\, C_{2}\in \xi$, then
     $\cap_{e\in EC_{1}} \G_{e}$ is conjugate in $\G$ to $\cap_{e\in
     EC_{2}} \G_{e}$, and we denote by $\cai(\xi)$ their common
     cardinality.

     Let $C\in\Re_{m}$ be such that $[C]_{\G}=\xi$.  By choosing each
     time a different starting edge, we obtain
     $\ell:=\ell(C)\equiv\ell(\xi)$ closed paths from $C$.  Denote
     them by $D_{1},\ldots,D_{\ell}$, and observe that any two of them
     can be $\G$-equivalent, $i.e.$ $D_{i}=\g(D_{j})$, for some
     $\g\in\G$, if and only if $\g\in \G_{C}$.  Moreover, if
     $\g\in\cap_{e\in EC} \G_{e}\subset \G_{C}$, then $\g(D_{i})
     =D_{i}$, for $i=1,\ldots,\ell$.  Therefore, there are only
     $k\equiv k(\xi):=\frac{\ell(\xi)\cai(\xi)}{\cs(\xi)}$ distinct
     $\G$-classes of closed paths generated by the $D_{i}$'s, and we
     denote them by $\pi_{1},\ldots,\pi_{k}$.

     Let $\pi$ be one of them, and observe that, for any
     $e\in\cf_{1}$, there are either no closed paths $D$ representing
     $\pi$ and such that $e_{0}(D)=e$, or there are
     $\frac{|\G_{e}|}{\cai(\xi)}$ distinct closed paths $D$
     representing $\pi$ and such that $e_{0}(D)=e$.  Indeed, if there is a
     closed path $D$ representing $\pi$ and such that $e_{0}(D)=e$, then
     any $\g\in \G_{e}$ generates a closed path $\g(D)$ representing
     $\pi$ and such that $e_{0}(\g(D))=e$, but, if $\g\in \cap_{e\in ED}
     \G_{e}$, then $\g(D)=D$. Hence, the claim is established.

     Let us now introduce a discrete measure on $\cf_{1}$.  Let us say
     that a $\G$-class of closed paths $\pi$ starts at $e\in\cf_{1}$
     if there is $D\in\pi$ such that $e_{0}(D)=e$.  Let us set, for
     $e\in\cf_{1}$, $\m_{\xi}(e)=1$, if $e$ is visited by some
     $\pi_{i}$, $i=1,\ldots,k$, and $\m_{\xi}(e)=0$, otherwise.  It is
     easy to see that $\m_{\xi}$ depends only on $\xi$ and is in
     particular independent of the representative $C$.  Observe that
     $\m_{\xi}(\cf_{1}) = k(\xi)$.

     Therefore, for any $e\in\cf_{1}$, we get
     $$
     |\set{D\in\Ci_{m}: [D]_{o,\G}=\xi, e_{0}(D)=e}| =
     \frac{\m_{\xi}(e)\cdot|\G_{e}|}{\cai(\xi)},
     $$
     and, finally,
     $$
     \sum_{e\in\cf_{1}} \frac{1}{|\G_{e}|}
     |\set{D\in\Ci_{m}: [D]_{o,\G}=\xi, e_{0}(D)=e}| =
     \frac{1}{\cai(\xi)}
     \sum_{e\in\cf_{1}} \m_{\xi}(e) = \frac{k(\xi)}{\cai(\xi)} =
     \frac{\ell(\xi)}{\cs(\xi)}.
     $$
 \end{proof}

 Define, for $\o\in \ell^{2}(EX)$, $e\in EX$,
 $$
 (T\o)(e) =\sum_{\begin{smallmatrix} t(e')=o(e) \\ e'\neq\ov{e}
 \end{smallmatrix} } \o(e').
 $$
 Then, we have

 \begin{Prop} \label{Lem:T}
     \itm{i} $T\in\cn_{1}(X,\G)$, $\|T\|\leq d-1$,

     \itm{ii} for $m\in\bn$, $T^{m}e =
     \sum_{\begin{smallmatrix}(e,e_{1},\ldots,e_{m}) \\ \text{proper
     path} \end{smallmatrix}} e_{m}$, for $e\in EX$,

     \itm{iii} $Tr_{\G}(T^{m}) = N^{\G}_{m} :=
     \sum_{[C]_{\G}\in[\Re_{m}]_{\G}}
     \frac{\ell([C]_{\G})}{\cs([C]_{\G})}$, the number of
     $\G$-equivalence classes of reduced cycles of length $m$.  Here,
     $Tr_{\G}$ is the trace on $\cn_{1}(X,\G)$ introduced in
     (\ref{eq:TraceOnEdges}).
 \end{Prop}
 \begin{proof}
     $(i)$, $(ii)$ are easy to check.

     $(iii)$ Using Lemma \ref{counting}, we obtain


     \begin{align*}
     Tr_{\G}(T^{m}) &= \sum_{e\in \cf_{1}} \frac{1}{|\G_{e}|}
     T^{m}(\tilde{e}, \tilde{e}) \\
     & = \sum_{e\in \cf_{1}} \frac{1}{|\G_{e}|}
     \sum_{\begin{smallmatrix}(e,e_{1},\ldots,e_{m-1},e)
     \\
     \text{reduced path} \end{smallmatrix}} 1 \\
     & = \sum_{e\in \cf_{1}} \frac{1}{|\G_{e}|}
     |\set{C\in\Ci_{m}: e_{0}(C)=e}| \\
     & = \sum_{[C]_\G \in [\Re]_\G} \sum_{e\in \cf_{1}}
     \frac{1}{|\G_{e}|} |\set{D\in\Ci_{m}: [D]_{0}\sim_\G C,
     e_{0}(D)=e}| \\
     &= N^{\G}_m.
     \end{align*}
 \end{proof}

\section{The Zeta function} \label{sec:Zeta}

 Before introducing the zeta function of an infinite periodic graph,
 we recall its definition for a finite $(q+1)$-regular graph $X$
 ($i.e.$ such that $\deg(v) =q+1$, for all $v\in VX$).  In that case,
 the Ihara zeta function $Z_X$ is defined by an Euler product of the
 form
 \begin{equation}
     Z_X(u) := \prod_{C\in \Pr} (1-u^{|C|})^{-1}, \ \text{ for }
     |u|<\frac{1}{q},
 \end{equation}
 where $\Pr$ is the set of prime cycles of $X$.  By way of comparison,
 recall that the Riemann zeta function is given by the Euler product
 \begin{equation}
     \zeta(s) := \prod_{p} (1-p^{-s})^{-1}, \ \text{ for } Re\ s>1,
 \end{equation}
 where $p$ ranges over all the rational primes.  To see the
 correspondence between $Z_X$ and $\zeta$, simply let $u:= q^{-s}$ and
 observe that $u^{|C|} = (q^{|C|})^{-s}$.  Also note that
 $|u|<\frac{1}{q}$ if and only if $Re\ s>1$.

 Let us now return to the case of periodic graphs and introduce the
 Ihara zeta function via its Euler product as well as show that this
 defines a holomorphic function in a suitable disc.

 \begin{Dfn}[Zeta function]
     Let $Z(u)=Z_{X,\G}(u)$ be given by
     $$
     Z_{X,\G}(u) := \prod_{[C]_{\G}\in [\Pr]_{\G}} (1-u^{|C|})^{
     -\frac{1}{ |\G_{C}| } },
     $$
     for $u\in\bc$ sufficiently small so that the infinite product
     converges.
 \end{Dfn}

 In the following proposition we let
 $$
 \Det_{\G}(B) := \exp\,\circ\,
 Tr_{\G}\,\circ\log(B), \quad \text{ for } B\in\cn_{1}(X,\G).
 $$
  We refer to Section
 \ref{sec:AnalyticDet} for more details. Formula $(iv)$ in the following theorem was first
 established in \cite{ClMS1}, although with a different proof.

 \begin{Thm}\label{Prop:power.series}
     \itm{i} $Z(u):=\prod_{[C]_{\G} \in [\Pr]_{\G}} (1-u^{|C|})^{
     -\frac{1}{ |\G_{C}| } }$ defines a holomorphic function in the
     open disc $\{u\in\bc: |u|<\frac{1}{d-1}\}$.

     \itm{ii} $u\frac{Z'(u)}{Z(u)} = \sum_{m=1}^{\infty}
     N^{\G}_{m}u^{m}$, for  $|u|<\frac{1}{d-1}$.

     \itm{iii} $Z(u) = \exp\left( \sum_{m=1}^{\infty}
     \frac{N^{\G}_{m}}{m}u^{m} \right)$, for  $|u|<\frac{1}{d-1}$.

     \itm{iv} $Z(u) = \det_{\G}(I-uT)^{-1}$, for  $|u|<\frac{1}{d-1}$.
 \end{Thm}
 \begin{proof}
     Observe that it follows from Proposition \ref{Lem:T} that
     $\sum_{m=1}^{\infty} \frac{N^{\G}}{m}u^{m}$ defines a function
     which is holomorphic in $\{u\in\bc: |u|<\frac{1}{d-1}\}$.  Moreover, for
     any $u\in\bc$ such that $|u|<\frac{1}{d-1}$,
     \begin{align*}
     \sum_{m=1}^{\infty} N^{\G}_{m} u^{m} & = \sum_{[C]_{\G}\in
     [\Re]_{\G}} \frac{\ell([C]_{\G})}{\cs([C]_{\G})}\, u^{|C|}
     \\
     & = \sum_{[C]_{\G}\in [\Pr]_{\G}} \sum_{m=1}^{\infty}
     \frac{|C|}{|\G_{C}|}\, u^{|C^{m}|} \\
     &= \sum_{[C]_{\G}\in
     [\Pr]_{\G}}\frac{1}{|\G_{C}|}\, \sum_{m=1}^{\infty} |C| u^{|C|m} \\
     &= \sum_{[C]_{\G}\in [\Pr]_{\G}} \frac{1}{|\G_{C}|}\,  u\frac{d}{du}
     \sum_{m=1}^{\infty} \frac{u^{|C|m}}{m} \\
     &= -\sum_{[C]_{\G}\in [\Pr]_{\G}} \frac{1}{|\G_{C}|}\, u\frac{d}{du}
     \log(1-u^{|C|}) \\
     & = u\frac{d}{du} \log Z(u),
     \end{align*}
     where, in the last equality, we have used uniform convergence on
     compact subsets of $\set{u\in\bc: |u|<\frac{1}{d-1}}$.  From what
     has already been proved, $(i)-(iii)$ follow.  Finally, for
     $|u|<\frac{1}{d-1}$, we have
    \begin{align*}
    \log Z(u) &= \sum_{m=1}^{\infty} \frac{N^{\G}_{m}}{m} u^{m} \\
    &= \sum_{m=1}^{\infty} \frac{1}{m} Tr_{\G}((Tu)^{m}) \\
    &= Tr_{\G}\left(\sum_{m=1}^{\infty} \frac{(Tu)^{m}}{m} \right) \\
    &= Tr_{\G}(-\log(I-uT)).
    \end{align*}
 \end{proof}

 \begin{exmp}
     Some examples of cycles with different stabilizers are shown in
     figures \ref{fig:cycle1}, \ref{fig:cycle2}.  They refer to the
     graph in figure \ref{fig:Example} which is the standard lattice
     graph $X=\bz^{2}$ endowed with the action of the group $\G$
      generated by the reflection along the $x$-axis and the
     translations by elements $(m,n)\in\bz^{2}$, acting as
     $(m,n)(v_{1},v_{2}):= (v_{1}+4m,v_{2}+4n)$, for
     $v=(v_{1},v_{2})\in VX=\bz^{2}$.
     \begin{figure}[ht]
     \centering \psfig{file=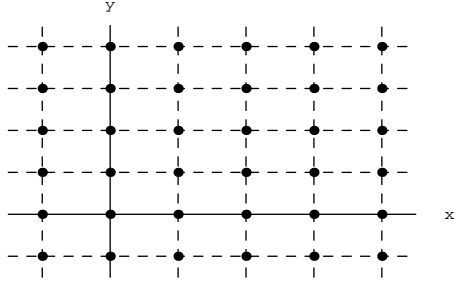,height=1.5in} \caption{A
     periodic graph} \label{fig:Example}
     \end{figure}
     \begin{figure}[ht]
     \centering
     \psfig{file=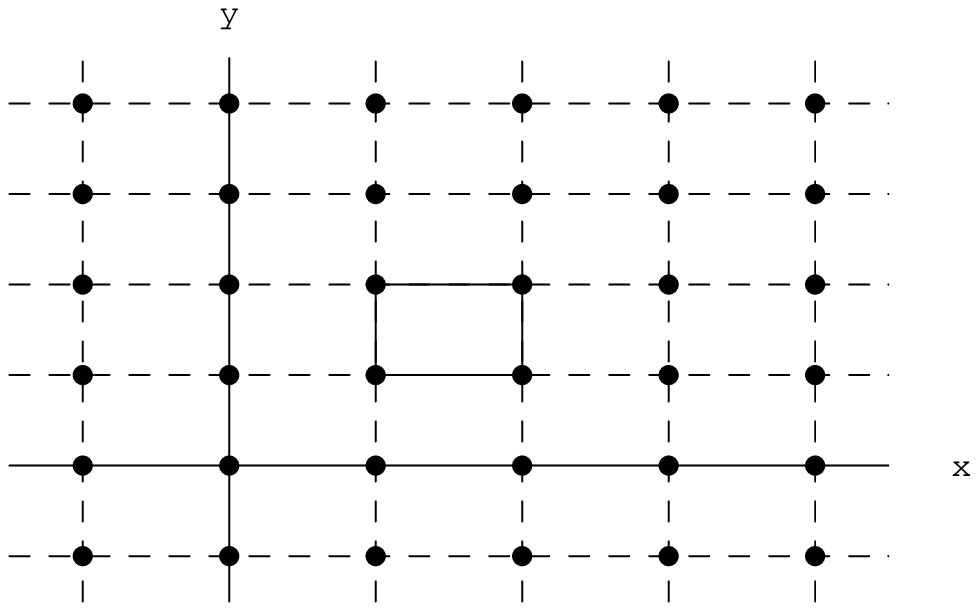,height=1.5in}
     \caption{A cycle with $|\G_{C}|=1$}
     \label{fig:cycle1}
     \end{figure}
     \begin{figure}[ht]
     \centering
     \psfig{file=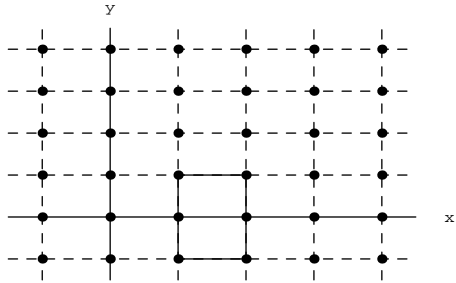,height=1.5in}
     \caption{A cycle with $|\G_{C}|=2$}
     \label{fig:cycle2}
     \end{figure}
 \end{exmp}

\section[Analytic determinant for von Neumann algebras]{An analytic
determinant for von Neumann algebras with a finite
trace}\label{sec:AnalyticDet}

 In this section, we define a determinant for a suitable class of not
 necessarily normal operators in a von Neumann algebra with a finite
 trace.  The results obtained are used in Section \ref{sec:DetFormula}
 to prove a determinant formula for the zeta function.

 In a celebrated paper \cite{FuKa}, Fuglede and Kadison defined a
 positive-valued determinant for finite factors ($i.e.$ von Neumann
 algebras with trivial center and finite trace).  Such a determinant is
 defined on all invertible elements and enjoys the main properties of
 a determinant function, but it is positive-valued.  Indeed, for an
 invertible operator $A$ with polar decomposition $A=UH$, where $U$ is a
 unitary operator and $H:= \sqrt{A^{*}A}$ is a positive self-adjoint
 operator, the Fuglede--Kadison determinant is defined by
 $$
 Det(A)=\exp\, \circ\ \tau\circ\log H,
 $$
 where $\log H$ may be defined via the functional calculus. Note, however, that the original definition was only given for a normalized trace.

 For the purposes of the present paper, we need a determinant which is
 an analytic function.  As we shall see, this can be achieved, but
 corresponds to a restriction of the domain of the determinant
 function and implies the loss of some important properties.  In particular,
 the product formula of the Fuglede--Kadison determinant only holds
 under certain restrictions in our case; see Propositions \ref{properties},
 \ref{prop:detProd}, \ref{prop:detAut} and \ref{prop:Binet}.

 Let
 $(\ca,\tau)$ be a von Neumann algebra endowed with a finite trace.
 Then, a natural way to obtain an analytic function is to define, for
 $A\in\ca$, $\Det_\t(A)=\exp\, \circ\ \tau\circ\log A$, where
 $$
 \log(A) := \frac{1}{2\pi i} \int_\Gamma \log \lambda (\lambda-A)^{-1}
 d\lambda,
 $$
 and $\Gamma$ is the boundary of a connected, simply connected region
 $\Omega$ containing the spectrum of $A$.  Clearly, once the branch of
 the logarithm is chosen, the integral above does not depend on
 $\Gamma$, provided $\G$ is given as above.

 Then a na\"{\i}ve way of defining $det$ is to allow all elements $A$
 for which there exists an $\Omega$ as above, and a branch of the
 logarithm whose domain contains $\Omega$.  Indeed, the following
 holds.

\begin{Lemma}
    Let $A$, $\Omega$, $\Gamma$ be as above, and $\varphi$, $\psi$ two
    branches of the logarithm such that both domains contain $\Omega$.
    Then
    $$
    \exp\, \circ\ \tau\circ\varphi(A) = \exp\, \circ\
\tau\circ\psi(A).
    $$
\end{Lemma}

\begin{proof}
    The function $\varphi(\lambda)-\psi(\lambda)$ is continuous and
    everywhere defined on $\Gamma$.  Since it takes its values in
    $2\pi i \mathbb{Z}$, it should be constant on $\Gamma$.  Therefore,
    \begin{align*}
     \exp\, \circ\ \tau\circ\varphi (A) & = \exp\, \circ\
    \tau\left(\frac{1}{2\pi i} \int_\Gamma 2\pi i n_{0}
    (\lambda-A)^{-1} d\lambda \right) \exp\, \circ\ \tau\circ\psi(A)\\
    &=\exp\, \circ\ \tau\circ\psi(A).
    \end{align*}
\end{proof}

 The problem with the previous definition is its dependence on the
 choice of $\Omega$.  Indeed, it is easy to see that when
 $A=\begin{pmatrix}1&0\\0&i\end{pmatrix}$ and we choose $\Omega$
 containing $\{e^{i\vartheta},\vartheta\in[0,\pi/2]\}$ and any
 suitable branch of the logarithm, we get $det(A)=e^{i\pi/4}$, if we
 use the normalized trace on $2\times 2$ matrices.  By contrast, if we
 choose $\Omega$ containing
 $\{e^{i\vartheta},\vartheta\in[\pi/2,2\pi]\}$ and a corresponding
 branch of the logarithm, we get $det(A)=e^{5i\pi/4}$.  Therefore, we
 make the following choice.

\begin{Dfn}
    Let $(\ca,\tau)$ be a von Neumann algebra endowed with a finite
    trace, and consider the subset $\ca_{0}=\{A\in\ca : 0\not\in
    \conv\sigma(A)\}$, where $\s(A)$ denotes the spectrum of $A$.  For
    any $A\in\ca_{0}$ we set
    $$
    \Det_\t(A)=\exp\, \circ\ \tau\circ\left(\frac{1}{2\pi i}
    \int_\Gamma \log \lambda (\lambda-A)^{-1} d\lambda\right),
    $$
    where $\Gamma$ is the boundary of a connected, simply connected
    region $\Omega$ containing $\conv\sigma(A)$, and $\log$ is a
    branch of the logarithm whose domain contains $\Omega$.
\end{Dfn}

\begin{Cor}\label{cor:det.analytic}
    The determinant function defined above is well defined and
    analytic on $\ca_{0}$.
\end{Cor}

We collect several properties of our determinant in the following
result.

\begin{Prop}\label{properties}
    Let $(\ca,\tau)$ be a von Neumann algebra endowed with a finite
    trace, and let $A\in\ca_{0}$.  Then

    \item[$(i)$] $\Det_\t(zA)=z^{\t(I)}\Det_\t(A)$, for any
    $z\in\mathbb{C}\setminus\{0\}$,

    \item[$(ii)$] if $A$ is normal, and $A=UH$ is its polar
    decomposition,
    $$\Det_\t (A)=\Det_\t(U)\Det_\t(H),$$

    \item[$(iii)$] if $A$ is positive, $\Det_\t(A)=Det(A)$, where the
    latter is the Fuglede--Kadison determinant.
\end{Prop}

\begin{proof}
    $(i)$ If, for a given $\vartheta_0\in[0,2\pi)$, the half-line $\{\rho e^{i\vartheta_0}\in\mathbb{C} :
    \r>0\}$ does not intersect $\conv\sigma(A)$, then the half-line
    $\{\rho e^{i(\vartheta_0+t)}\in\mathbb{C} : \r>0\}$ does not
    intersect $\conv\sigma(zA)$, where $z=re^{it}$.  If $\log$ is the
    branch of the logarithm defined on the complement of the real
    negative half-line, then $\varphi(x)=i(\vartheta_{0}-\pi) +
    \log(e^{-i(\vartheta_{0}-\pi)}x)$ is suitable for defining
    $\Det_\t(A)$, while $\psi(x)=i(\vartheta_{0}+t-\pi) +
    \log(e^{-i(\vartheta_{0}+t-\pi)}x)$ is suitable for defining
    $\Det_\t(zA)$.  Moreover, if $\Gamma$ is the boundary of a
    connected, simply connected region $\Omega$ containing
    $\conv\sigma(A)$, then $z\Gamma$ is the boundary of a connected,
    simply connected region $z\Omega$ containing $\conv\sigma(zA)$.
    Therefore,
    \begin{align*}
    \Det_\t(zA) &= \exp\, \circ\ \tau\left(\frac{1}{2\pi i}
    \int_{z\Gamma} \psi(\lambda) (\lambda-zA)^{-1}
    d\lambda\right)\\
    &= \exp\, \circ\ \tau\left(\frac{1}{2\pi i} \int_{\Gamma}
    (i(\vartheta_{0}+t-\pi) + \log(e^{-i(\vartheta_{0}+t-\pi)}
    re^{it}\mu)) (\mu-A)^{-1} d\mu\right)\\
    &= \exp\, \circ\ \tau\left((\log r + it)I+\frac{1}{2\pi i}
    \int_{\Gamma} \varphi(\mu) (\mu-A)^{-1} d\mu\right)\\
    &= z^{\t(I)} \Det_\t(A).
    \end{align*}
    $(ii)$ When $A=UH$ is normal, $U=\int_{[0,2\pi]} e^{i\vartheta}\
    du(\th)$, $H=\int_{[0,\infty)}r\ dh(r)$, then $ A =
    \int_{[0,\infty)\times[0,2\pi]} r e^{i\vartheta} \ d(h(r)\otimes
    u(\th))$.  The property $0\not\in\conv\sigma(A)$ is equivalent to
    the fact that the support of the measure $d(h(r)\otimes u(\th))$
    is compactly contained in some open half-plane $$\{\rho
    e^{i\vartheta} : \rho>0, \vartheta \in (\vartheta_{0} - \pi/2,
    \vartheta_{0} +\pi/2)\},$$ or, equivalently, that the support of
    the measure $dh(r)$ is compactly contained in $(0,\infty)$, and
    the support of the measure $d u(\th)$ is compactly contained in
    $(\vartheta_{0} - \pi/2, \vartheta_{0} +\pi/2)$.  Therefore,
    $A\in\ca_{0}$ is equivalent to $U,H\in\ca_{0}$.  Then $$\log A =
    \int_{[0,\infty) \times (\vartheta_{0} - \pi/2, \vartheta_{0}
    +\pi/2)} (\log r + i\vartheta) \ d(h(r)\otimes u(\th)),$$ which
    implies that
    \begin{align*}
    \Det_\t(A) &= \exp\, \circ\ \tau\left(\int_{0}^{\infty} \log
    r\ dh(r) + \int_{\vartheta_{0} - \pi/2}^{\vartheta_{0} +\pi/2}
    i\vartheta \ du(\th)\right) \\
    &= \Det_\t(U)\cdot \Det_\t(H).
    \end{align*}
    $(iii)$  This follows by the  argument given in $(ii)$.
\end{proof}

\begin{rem}
    We note that the above defined determinant function strongly
    violates the product property $\Det_\t(AB)=\Det_\t(A)\Det_\t(B)$.
    Indeed, the fact that $A,B\in\ca_{0}$ does not imply
    $AB\in\ca_{0}$, as is seen e.g. by taking
    $A=B=\begin{pmatrix}1&0\\0&i\end{pmatrix}$.  Moreover, even if
    $A,B,AB\in\ca_{0}$ and $A$ and $B$ commute, the product property
    may be violated, as is shown by choosing
    $A=B=\begin{pmatrix}1&0\\0&e^{3i\pi/4}\end{pmatrix}$, and using
    the normalized trace on $2\times 2$ matrices.
\end{rem}

 \begin{Prop} \label{prop:detProd}
     Let $(\ca,\tau)$ be a von Neumann algebra endowed with a finite
     trace, and let $A,B\in\ca$.  Then, for sufficiently small
     $u\in\bc$, we have
     $$
     \Det_{\t}((I+uA)(I+uB)) = \Det_{\t}(I+uA) \Det_{\t}(I+uB).
     $$
 \end{Prop}
 \begin{proof}
     The proof is inspired by that of Lemma 3 in \cite{FuKa}.  Let us
     write $a:= \log(I+uA)$, $b:= \log(I+uB)\in\ca$, and let $c(t) :=
     e^{ta}e^b$, $t\in[0,1]$.  As $\|a\| \leq -\log(1-|u|\|A\|)$, and
     $\|b\| \leq -\log(1-|u|\|B\|)$, we get
     \begin{align*}
     \| c(t)-1\| & = \|e^{ta}-e^{-b}\| \|e^b\|\\
     &\leq e^{\|b\|} \left( e^{\|a\|} + e^{\|b\|} -2\right)\\
     & \leq \frac{1}{1-|u|\|B\|} \left( \frac{1}{1-|u|\|A\|} +
     \frac{1}{1-|u|\|B\|} -2\right) <1,
     \end{align*}
     for all $t\in[0,1]$, if we choose $|u|$ sufficiently small;
     hence, $c(t)\in\ca_0$ for all $t\in[0,1]$.  Now apply Lemma 2 in
     \cite{FuKa} which gives
     $$
     \t( \frac{d}{dt} \log c(t) ) = \t(c(t)^{-1}c'(t)) =
     \t(e^{-b}e^{-ta}ae^{ta}e^b) = \t(a).
     $$
     Therefore, after integration for $t\in[0,1]$, we obtain $\t(\log
     c(1)) -\t(\log c(0)) = \t(a)$, which means
     \begin{align*}
     \t\bigl(\log((I+uA)(I+uB))\bigr) &= \t(\log c(1)) = \t(a)+\t(b) \\
     &= \t\bigl(\log(I+uA)\bigr)+ \t\bigl(\log(I+uB)\bigr),
     \end{align*}
     and hence implies the claim.
 \end{proof}

 \begin{Prop} \label{prop:detAut}
     Let $(\ca,\tau)$ be a von Neumann algebra endowed with a finite
     trace.  Further, let $A\in\ca$ have a bounded inverse, and let
     $T\in\ca_{0}$.  Then
     $$
     \Det_{\t}(ATA^{-1}) = \Det_{\t} T.
     $$
 \end{Prop}
 \begin{proof}
     Indeed, for any polynomial $p$, we have $p(ATA^{-1}) =
     Ap(T)A^{-1}$.  Applying the Stone--Weierstrass theorem on the
     compact set $\s(ATA^{-1})=\s(T)$, we obtain $\log(ATA^{-1}) =
     A\log(T)A^{-1}$, from which the result follows.
 \end{proof}

  \begin{Prop} \label{prop:Binet}
     Let $(\ca,\tau)$ be a von Neumann algebra endowed with a finite
    trace, and let $\displaystyle{ T= \left(
    \begin{array}{cc}
        T_{11} & T_{12}  \\
        0 & T_{22}
    \end{array} \right) \in Mat_{2}(\ca) }$, with
    $T_{ii}\in\ca$ such that $\s(T_{ii})\subset B_1(1):= \set{z\in\bc : |z-1|<1}$,
    for $i=1,2$.  Then
     $$
     \Det_{\t}(T) = \Det_{\t}(T_{11}) \Det_{\t}(T_{22}).
     $$
 \end{Prop}
 \begin{proof}
     Indeed, for any $k\in\bn\cup\set{0}$,
     $$
     T^{k} = \left(
    \begin{array}{cc}
        T_{11}^{k} & B_{k}  \\
        0 & T_{22}^{k}
    \end{array} \right),
     $$
     for some $B_{k}\in\ca$, so that, for any polynomial $p$,
     $$
     p(T) = \left(
    \begin{array}{cc}
        p(T_{11}) & B  \\
        0 & p(T_{22})
    \end{array} \right),
     $$
     for some $B\in\ca$.  It is easy to see that $\s(T)\subset\s(T_{11})\cup \s(T_{22})\subset B_1(1)$. Hence, applying the Stone--Weierstrass theorem on
     the compact set $\s(T)$, we obtain
     $$
     \log(T) = \left(
    \begin{array}{cc}
        \log(T_{11}) & C  \\
        0 & \log(T_{22})
    \end{array} \right),
     $$
     for some $C\in\ca$.  Therefore,
     $$
     \Det_{\t}(T) = \exp\, \circ \t\ \circ\, \log(T) = \exp\bigl(
     \t(\log(T_{11})) + \t(\log(T_{22}))\bigr) = \Det_{\t}(T_{11})
     \Det_{\t}(T_{22}),
     $$
     as desired.
 \end{proof}

 \begin{Cor} \label{cor:Binet}
     Let $\G$ be a discrete group, $\pi_{1},\, \pi_{2}$ unitary
     representations of $\G$, and $\t_{1},\,\t_{2}$ finite traces on
     $\pi_{1}(\G)'$ and $\pi_{2}(\G)'$, respectively.  Let $\pi:=
     \pi_{1}\oplus \pi_{2}$, $\t:=\t_{1}+\t_{2}$, $\displaystyle{ T=
     \left(
     \begin{array}{cc}
        T_{11} & T_{12}  \\
     0 & T_{22} \end{array} \right) \in \pi(\G)' }$, with $\s(T_{ii})\subset B_1(1)= \set{z\in\bc : |z-1|<1}$, for $i=1,2$. Then
     $$
     \Det_{\t}(T) = \Det_{\t_{1}}(T_{11}) \Det_{\t_{2}}(T_{22}).
     $$
 \end{Cor}
 \begin{proof}
     It is similar to the proof of Proposition \ref{prop:Binet}.
 \end{proof}

\section{The determinant formula} \label{sec:DetFormula}

 In this section, we prove the main result in the theory of the Ihara
 zeta functions, which says that $Z$ is the reciprocal of a
 holomorphic function, which, up to a factor, is the determinant of a
 deformed Laplacian on the graph.  We first need some technical
 results.

 Let us denote by $A$ the adjacency matrix of $X$, $i.e.$
 $\displaystyle{(Af)(v) = \sum_{w\sim v}f(w)}$, $f\in\ell^{2}(VX)$.
 Then (by \cite{Mohar}, \cite{MoWo}) $\|A\|\leq d:=\sup_{v\in VX}
 \deg(v) <\infty$, and it is easy to see that $A\in\cn_{0}(X,\G)$.
 Introduce $(Qf)(v) := (\deg(v)-1)f(v)$, $v\in VX$,
 $f\in\ell^{2}(VX)$, and $\D(u):= I-uA+u^{2}Q\in\cn_{0}(X,\G)$, for
 $u\in\bc$.  Let us recall that $d:=\sup_{v\in VX} \deg(v)$, and set
 $\a:= \frac{d+\sqrt{d^{2}+4d}}{2}$.  Then

 \begin{Thm}[Determinant formula] \label{Thm:det.formula}
     We have
     $$
     Z_{X,\G}(u)^{-1} = (1-u^{2})^{-\chi^{(2)}(X)}
     \Det_{\G}(\D(u)),\ \text{ for }
     |u|<\frac{1}{\a},
     $$
     where $\displaystyle \chi^{(2)}(X) := \sum_{v\in\cf_0}
     \frac{1}{|\G_{v}|} - \frac12 \sum_{e\in \cf_1}
     \frac{1}{|\G_{e}|}$ is the $L^{2}$-Euler
     characteristic of $(X,\G)$, as introduced in {\rm \cite{ChGr}}.
 \end{Thm}

 This theorem was first proved in \cite{ClMS1} and is based on formula $(iv)$
 in Theorem \ref{Prop:power.series} and the equality
 $\Det_\G(I-uT) = (1-u^{2})^{-\chi^{(2)}(X)} \Det_{\G}(\D(u))$, for $ |u|<\frac{1}{\a}$.
 The main difference with their proof is that we use an analytic determinant and operator-valued
 analytic functions instead of Bass' noncommutative determinant \cite{Bass} and formal power series of operators.

 We first prove two lemmas. Define, for $f\in\ell^{2}(VX),\
 \o\in\ell^{2}(EX)$,
 \begin{align*}
     (\de_{0}f)(e) &:= f(o(e)),\ e\in EX\\
     (\de_{1}f)(e) &:= f(t(e)),\ e\in EX\\
     (\s\o)(v) &:= \sum_{o(e)=v} \o(e),\ v\in VX\\
     (J\o)(e)&:= \o(\ov{e}),\ e\in EX,
 \end{align*}
 and use the short-hand notation $I_{V} := Id_{\ell^{2}(VX)}$ and
 $I_{E}:= Id_{\ell^{2}(EX)}$.

 \begin{Lemma}
     \itm{i} $J\de_{1}=\de_{0}$,

     \itm{ii} $\s\l_{1}(\g)=\l_{0}(\g)\s$,
     $\de_{i}\l_{0}(\g)=\l_{1}(\g)\de_{i}$, $i=0,1$, $\g\in\G$,

     \itm{iii} $\s\de_{0} = I+Q$,

     \itm{iv} $\s\de_{1} = A$,

     \itm{v} $\de_{0}\s = JT+I_{E}$,

     \itm{vi} $\de_{1}\s = T+J$,

     \itm{vii} $(I_{E}-uJ)(I_{E}-uT) = (1-u^{2})I_{E} -
     u\de_{1}\s+u^{2}\de_{0}\s$.
 \end{Lemma}
 \begin{proof}
     Let $f\in\ell^2(VX)$, $v\in VX$.  Then
     \begin{align*}
     (\s\de_0f)(v) & = \sum_{o(e)=v} (\de_0f)(e) = \sum_{o(e)=v}
     f(o(e)) = (1+Q(v,v))f(v) \\
     (\s\de_1f)(v) & = \sum_{o(e)=v} (\de_1f)(e) = \sum_{o(e)=v}
     f(t(e)) =(Af)(v).
     \end{align*}
     Moreover, for $\o\in\ell^2(EX)$, $e\in EX$, we have
     \begin{align*}
     (\de_1\s\o)(e) & = (\s\o)(t(e)) = \sum_{o(e')=t(e)}
     \o(e') = (T\o)(e) + (J\o)(e)\\
     \de_0\s & = J\de_1\s  = J(T+J) = JT+I_E.
     \end{align*}
 The rest of the proof is clear.
 \end{proof}

 \medskip
 Let us now consider the direct sum of the unitary representations $\l_0$ and $\l_1$,  namely
 $\l(\g) := \l_0(\g)\oplus \l_1(\g) \in \cb(\ell^2(VX)\oplus
 \ell^2(EX))$.  Then, the von Neumann algebra $\l(\G)' := \set{S\in
 \cb(\ell^2(VX)\oplus \ell^2(EX)) : S\l(\g) = \l(\g)S, \ \g\in\G}$
 consists of operators $S= \left(\begin{array}{cc} S_{00} & S_{01} \\
 S_{10} & S_{11} \\
 \end{array}\right)$, where $S_{ij}\l_j(\g) = \l_i(\g)S_{ij},\
 \g\in\G,\ i,j=0,1$, so that $S_{ii}\in\l_i(\G)'\equiv \cn_i(X,\G)$,
 $i=0,1$.  Hence $\l(\G)'$ inherits a trace given by
 \begin{equation}\label{eq:TraceOnDirectSum}
      Tr_\G \left( \begin{array}{cc}
                      S_{00} & S_{01} \\
                      S_{10} & S_{11} \\
                   \end{array}\right)
      := Tr_\G(S_{00}) + Tr_\G(S_{11}).
 \end{equation}

 Introduce
 $$
 \cl(u):=\left(
 \begin{array}{cc}
     (1-u^{2})I_{V} & 0  \\
     u\de_{0}-\de_{1} & I_{E}
 \end{array}\right) \text{ and }
 \cam(u) :=\left(
 \begin{array}{cc}
     I_{V} & u\s  \\
     u\de_{0}-\de_{1} & (1-u^{2})I_{E}
 \end{array}\right),
 $$
 which both belong to $\l(\G)'$. Then, we have

 \begin{Lemma}\label{lem:Technical}
     \itm{i} $\displaystyle{ \cam(u)\cl(u) = \left(
     \begin{array}{cc}
         \D(u) & u\s  \\
         0 & (1-u^{2})I_{E}
     \end{array}\right) }$,

     \itm{ii} $\displaystyle{ \cl(u)\cam(u) = \left(
     \begin{array}{cc}
         (1-u^{2})I_{V} & (1-u^{2})u\s  \\
         0 & (I_{E}-uJ)(I_{E}-uT)
     \end{array}\right) }$.

    \noindent Moreover, for $|u|$ sufficiently small,
     \itm{iii} $\cl(u)$, $\cam(u)$ are invertible, with a bounded
     inverse,

     \itm{iv} $\displaystyle{ \Det_{\G}(\cam(u)\cl(u)) =
     (1-u^{2})^{Tr_{\G}(I_{E})}\Det_{\G}(\D(u)) }$,

     \itm{v} $\displaystyle{ \Det_{\G}(\cl(u)\cam(u)) =
     (1-u^{2})^{Tr_{\G}(I_{V})-\frac12 Tr_{\G}(I_{E})}
     \Det_{\G}(I_{E}-uT) }$.
 \end{Lemma}
 \begin{proof}
 The formulas for $\cam(u)\cl(u)$ and $\cl(u)\cam(u)$ follow from the
 previous lemma.  Moreover, for $|u|$ sufficiently small,
 $\s( \D(u))$, $\s((1-u^{2})I_{E})$, $\s((1-u^{2})I_{V})$ and
 $\s((I_{E}-uJ)(I_{E}-uT))\subset B_1(1)= \set{z\in\bc : |z-1|<1}$,
 hence $\s( \cam(u)\cl(u) )$ and $\s(\cl(u)\cam(u))\subset B_1(1)$,
 as in the proof of Proposition \ref{prop:Binet}. Therefore,
 $\cl(u)$ and $\cam(u)$ are invertible, with a bounded inverse, for $|u|$ sufficiently
 small.   By Propositions
 \ref{properties} $(i)$, \ref{prop:detProd} and Corollary \ref{cor:Binet}, we
 obtain
 \begin{align*}
    \Det_{\G}(\cam(u)\cl(u)) & =
    \Det_{\G}(\D(u))\Det_{\G}((1-u^{2})I_{E}) \\
    & = (1-u^{2})^{Tr_{\G}(I_{E})} \Det_{\G}(\D(u))
 \end{align*}
 and
 \begin{align*}
    \Det_{\G}(\cl(u)\cam(u)) & =
    \Det_{\G}((1-u^{2})I_{V})\Det_{\G}(I_{E}-uJ)\Det_{\G}(I_{E}-uT) \\
    & = (1-u^{2})^{Tr_{\G}(I_{V})}
    \Det_{\G}(I_{E}-uJ)\Det_{\G}(I_{E}-uT).
 \end{align*}
 Moreover, we have $\Det_{\G}(I_{E}-uJ) = (1-u^{2})^{\frac{1}{2}\,
 Tr_{\G}(I_{E})}$.  Indeed, using $J$ to identify $\ell^{2}(E^{-}X)$
 with $\ell^{2}(E^{+}X)$, we obtain a representation $\r$ of
 $\cb(\ell^{2}(EX))$ onto $Mat_{2}(\cb(\ell^{2}(E^{+}X)))$, under
 which $\displaystyle{ \r(J)= \left(
 \begin{array}{cc}
     0 & I  \\
     I & 0
 \end{array}\right), \r(I_{E}) = \left(
 \begin{array}{cc}
     I & 0  \\
     0 & I
 \end{array}\right) }$. Hence, by Propositions \ref{prop:detProd} and
\ref{prop:Binet},
 \begin{align*}
     \Det_{\G}(I_{E}-uJ) & = \Det_{\G}(\r(I_{E}-uJ)) \\
     & = \Det_{\G} \left(
 \begin{array}{cc}
     I & -uI  \\
     -uI & I
 \end{array}\right) \\
 & = \Det_{\G} \left(\left(
 \begin{array}{cc}
     I & uI  \\
     0 & I
 \end{array}\right) \left(
 \begin{array}{cc}
     I & -uI  \\
     -uI & I
 \end{array}\right) \right)\\
 & = \Det_{\G} \left(
 \begin{array}{cc}
     (1-u^{2})I & 0  \\
     -uI & I
 \end{array}\right) \\
 & = (1-u^{2})^{Tr_{\G}(I)} \\
 & = (1-u^{2})^{\frac{1}{2}\, Tr_{\G}(I_{E})}.
 \end{align*}
 \end{proof}

 \noindent {\it Proof} (of Theorem \ref{Thm:det.formula}).

     Let us observe that, for sufficiently small $|u|$, we have
     $$\cam(u)\cl(u) = \cam(u) \cl(u)\cam(u) \cam(u)^{-1},$$ so that,
     by Proposition \ref{prop:detAut}, we get
     $\Det_{\G}(\cl(u)\cam(u)) = \Det_{\G}(\cam(u)\cl(u))$.
     Therefore, the claim follows from Lemma \ref{lem:Technical}
     $(iv)$ and $(v)$, equations (\ref{eq:TraceOnVertices}) and
     (\ref{eq:TraceOnEdges}) and Theorem \ref{Prop:power.series}.
 \endproof

\section{Functional equations}\label{sec:functEq}

 In this section, we obtain several functional equations for the Ihara
 zeta functions of $(q+1)$-regular graphs, $i.e.$ graphs with
 $\deg(v)=q+1$, for any $v\in VX$, on which $\G$ acts freely [$i.e.$
 $\G_{v}$ is trivial, for $v\in VX$] and with finite quotient [$i.e.$
 $B:=X/\G$ is a finite graph].  The various functional equations
 correspond to different ways of completing the zeta functions, as is done in \cite{StTe} for finite graphs. We extend here to non necessarily simple graphs the results  contained in \cite{GILa02}.

\begin{Lemma} \label{prop:holomorphy}
    Let $X$ be a $(q+1)$-regular graph, on which $\G$ acts freely and
    with finite quotient $B:= X/\G$.  Let $\D(u) := (1+qu^2)I-uA$. Then

    \itm{i} $\chi^{(2)}(X)=\chi(B)= |V(B)|(1-q)/2\in\bz$,

    \itm{ii} $\displaystyle Z_{X,\G}(u) = (1-u^2)^{\chi(B)}
    \Det_{\G}((1+qu^2)I-uA)^{-1}$, for $|u| < \frac{1}{q}$,

    \itm{iii} by using the determinant formula in $(ii)$, $Z_{X,\G}$ can
    be extended to a function holomorphic at least in the open set
    $$
    \O:=\br^2 \setminus \left(\set{(x,y)\in\br^2: x^2+y^2=\frac{1}{q}}
    \cup \set{(x,0)\in\br^2: \frac{1}{q}\leq |x|\leq 1 }\right).
    $$
    See figure \ref{fig:Omega}.
     \begin{figure}[ht]
     \centering
     \psfig{file=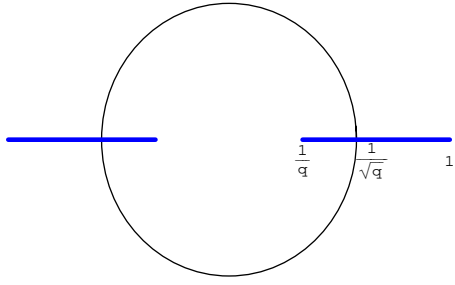,height=1.5in}
     \caption{The open set $\Omega$}
     \label{fig:Omega}
     \end{figure}

     \itm{iv} $\displaystyle \Det_\G\Bigl(\D ( \frac{1}{qu})\Bigr) = (qu^2)^{-|VB|} \Det_\G(\D(u))$, for $u\in\O\setminus \set{0}$.
\end{Lemma}
\begin{proof}
    $(i)$ This follows by a simple computation.

    $(ii)$ This follows from $(i)$.

    $(iii)$ Let us observe that
    $$
    \s(\D(u)) = \set{1+qu^2-u\l: \l\in\s(A)} \subset \set{1+qu^2-u\l:
    \l\in[-d,d]}.
    $$
    It follows that $0\not\in\conv\s(\D(u))$ at least
    for $u\in\bc$ such that $1+qu^2-u\l\neq0$ for $\l\in[-d,d]$, that
    is for $u=0$ or $\frac{1+qu^2}{u}\not\in[-d,d]$, or equivalently,
    at least for $u\in\O$.  The rest of the proof follows from
    Corollary \ref{cor:det.analytic}.

    $(iv)$ This follows from Proposition \ref{properties} $(i)$ and the fact that $Tr_\G(I_V) = |VB|$.
\end{proof}

\medskip

The question whether the extension of the domain of $Z_{X,\G}$ by means of the determinant formula is compatible with an analytic extension from the defining domain is a non-trivial issue, see  the recent paper by Clair \cite{Clair}.

 \begin{Thm} [Functional equations] \label{thm:FunctEq}
     Let $X$ be a $(q+1)$-regular graph, on which $\G$ acts freely and
     with finite quotient $B:= X/\G$.  Then, for all $u\in\O$, we have

     \itm{i} $\La_{X,\G}(u) := (1-u^{2})^{-\chi(B)}(1-u^{2})^{|VB|/2}
     (1-q^{2}u^{2})^{|VB|/2} Z_{X,\G}(u) =
     -\La_{X,\G}\Bigl(\frac{1}{qu}\Bigr)$,

     \itm{ii} $\xi_{X,\G}(u) := (1-u^{2})^{-\chi(B)} (1-u)^{|VB|} (1-qu)^{|VB|}
     Z_{X,\G}(u) = \xi_{X,\G}\Bigl(\frac{1}{qu}\Bigr)$,

     \itm{iii} $\Xi_{X,\G}(u) := (1-u^{2})^{-\chi(B)} (1+qu^{2})^{|VB|}
     Z_{X,\G}(u) = \Xi_{X,\G}\Bigl(\frac{1}{qu}\Bigr)$.
 \end{Thm}
 \begin{proof}
 They all follow from Lemma \ref{prop:holomorphy} $(iv)$ by a straightforward computation. We prove $(i)$ as an example.
     \begin{align*}
     \La_{X}(u) & = (1-u^{2})^{|VB|/2} (1-q^{2}u^{2})^{|VB|/2}
     \Det_{\G}(\D(u))^{-1} \\
     &= u^{|VB|}\Bigl(\frac{q^{2}}{q^{2}u^{2}}-1 \Bigr)^{|VB|/2} (qu)^{|VB|}
     \Bigl(\frac{1}{q^{2}u^{2}}-1 \Bigr)^{|VB|/2}
     \frac{1}{(qu^{2})^{|VB|}}\Det_{\G}\Bigl(
     \D(\frac{1}{qu}) \Bigr)^{-1}\\
     &= -\La_{X}\Bigl(\frac{1}{qu}\Bigr).
     \end{align*}
  \end{proof}

 \begin{rem}
     Recall that a key property of the Riemann zeta function $\zeta$
     is that its meromorphic continuation satisfies a functional
     equation $\xi(s)=\xi(1-s)$, for all $s\in\bc$, where $\xi(s) :=
     \pi^{-s/2} \G(s/2) \zeta(s)$ denotes the completion of $\zeta$
     and $\G$ is the usual Gamma function.  Likewise, in Theorem
     \ref{thm:FunctEq}, any of the functional equations relates the
     values of the corresponding completed Ihara zeta function at $s$
     and $1-s$, provided we set $u=q^{-s}$, as was explained at the
     beginning of Section \ref{sec:Zeta}.  Note that $\frac{1}{qu} =
     \frac{1}{q^{1-s}}$.
 \end{rem}

\section[Approximation for amenable graphs]{Approximation by finite
graphs in the amenable case}

 In this section, we show that the zeta function of a graph, endowed
 with a free and cofinite action of a discrete amenable group of
 automorphisms, is the limit of the zeta functions of a (suitable)
 sequence of finite subgraphs, thus answering in the affirmative a
 question raised by Grigorchuk and $\dot{\text{Z}}$uk in \cite{GrZu}.

 Before doing that, we establish a result which is considered folklore by
 specialists. Roughly speaking, it states that a $\G$-space is amenable
 if $\G$ is an amenable group, where a space is said to be amenable
 if it possesses a regular exhaustion. Such a result was stated by
 Cheeger and Gromov in \cite{ChGr} for CW-complexes and was proved by
 Adachi and Sunada in \cite{AdSu} for covering manifolds. We give here
 a proof in the case of covering graphs.

 Throughout this section, $X$ is a connected, countably infinite graph,
 and $\G$ is a countable discrete amenable group of automorphisms
 of $X$, which acts on $X$ freely [i.e., any $\g\neq id$ has no
 fixed-points], and cofinitely [i.e., $B:=X/\G$ is a finite graph].

 A fundamental domain for the action of $\G$ on $X$ can be constructed
 as follows.  Let $B=(VB,EB)$ be the quotient graph, and $p:X\to B$
 the covering map.  Let $EB = \set{e_{1},\ldots,e_{k}}$, where the
 edges have been ordered in such a way that, for each
 $i\in\set{1,\ldots,k}$, $e_{i}$ has at least a vertex in common with
 some $e_{j}$, with $j<i$.  Choose $\ee_{1}\in EX$ such that
 $p(\ee_{1})=e_{1}$.  Assume $\ee_{1},\ldots,\ee_{i}$ have already
 been chosen in such a way that $p(\ee_{j})=e_{j}$, for
 $j=1,\ldots,i$, and, for any such $j$, $\ee_{j}$ has at least a
 vertex in common with some $e_{h}$, with $h<j$.  Let $e_{i+1}\in EB$
 have a vertex in common with $e_{j}$, for some $j\in\set{1,\ldots,i}$
 and choose $\ee_{i+1}\in VX$ such that $p(\ee_{i+1})=e_{i+1}$ and
 $\ee_{i+1}$ has a vertex in common with $\ee_{j}$.  This completes
 the induction.  Let $EF := \set{\ee_{1},\ldots,\ee_{k}}$ and $VF :=
 \set{o(\ee_{1}),\ldots,o(\ee_{k})} \cup
 \set{t(\ee_{1}),\ldots,t(\ee_{k})}$, so that $F=(VF,EF)$ is a
 connected finite subgraph of $X$ which does not contain any
 $\G$-equivalent edges.  Then, $F$ is said to be a {\it fundamental
 domain} for the action of $\G$ on $X$.

 \begin{Dfn}\label{def:amenableGraph}
     Let $X$ be a countably infinite graph and $\G$ a countable
     discrete amenable group of automorphisms of $X$, which acts on
     $X$ freely and cofinitely; further, let $F$ be a corresponding
     fundamental domain.  A sequence $\set{K_n: n\in\bn}$ of finite
     subgraphs of $X$ is called an {\it amenable exhaustion} of $X$ if
     the following conditions hold:

     \itm{i} $K_{n} = \cup_{\g\in E_{n}} \g F$, where $E_{n}\subset
     \G$, for all $n\in\bn$,

     \itm{ii} $\cup_{n\in\bn} K_{n} =X$,

     \itm{iii} $K_{n}\subset K_{n+1}$, for all $n\in\bn$,

     \itm{iv} if $\cf K_{n}:= \set{v\in VK_{n} : d(v, VX\setminus
     VK_{n})=1}$, then $\displaystyle{ \lim_{n\to\infty} \frac{|\cf
     K_{n}|}{|VK_{n}|} =0}$.

     Then $X$ is called an {\it amenable graph} if it possesses an amenable
     exhaustion.
 \end{Dfn}

 \begin{Thm}\label{thm:amenable}
     Let $X$ be a connected, countably infinite graph, $\G$ be
     a countable discrete amenable subgroup of automorphisms of $X$
     which acts on $X$ freely and cofinitely and let $F$ be a
     corresponding fundamental domain.  Then $X$ is an amenable graph.
 \end{Thm}
 \begin{proof}
     The proof is an adaptation of a proof by Adachi and Sunada in
     the manifold case, see \cite{AdSu}.

 The finite set $A:= \set{\g\in\G: dist(\g F,F)\leq 1}$ is symmetric
 [$i.e.$ $\g\in A \iff \g^{-1}\in A$], generates $\G$ as a group, and
 contains the unit element.  Introduce the Cayley graph $\cc(\G,A)$,
 whose vertices are the elements of $\G$, and, by definition, there is
 one edge from $\g_{1}$ to $\g_{2}$ iff $\g_{1}^{-1}\g_{2}\in A$.  A
 subset $E\subset V\cc(\G,A)$ is said to be connected if, for any pair
 of distint vertices of $E$, there is a path in $\cc(\G,A)$, joining
 those two vertices, and consisting only of vertices of $E$.

 From \cite{Adachi}, Theorem 4, it follows that there is a sequence
 $\set{E_{j}}_{j\in\bn}$ of connected finite subsets of $\G$ such that
 \begin{eqnarray*}
     \cup_{j\in\bn} E_{j} =\G, \qquad E_{j}\subset E_{j+1}, \ \forall
     j\in\bn, \\
     \frac {|E_{j}\cdot A\setminus E_{j}|}{|E_{j}|} \leq
     \frac{1}{j|A|},\ \forall j\in\bn,
 \end{eqnarray*}
 where, for any $U_{1}$, $U_{2}\subset\G$, we set $U_{1}\cdot U_{2} =
 \set{\g_{1}\g_{2}: \g_{i}\in U_{i}, i=1,2}$.

 For each $n\in\bn$, let $K_{n}:= \cup_{\g\in E_{n}} \g F$.  Then
 $\set{K_{n}: n\in\bn}$ satisfies the claim.  Indeed, let $b:= |VF|$
 and $a:= |\cf_{0}|$, so that $a|E_{n}| \leq |VK_{n}| \leq b|E_{n}|$,
 $n\in\bn$.  Moreover, for any $n\in\bn$, we have
 $$
 \cf K_{n} \subset \cup_{\g\in U_{n}} \g F,
 $$
 where $U_{n} := \set{\g\in E_{n}: \text{ there is } \d\in A \text{
 such that } \g \d \not\in E_{n}}$.   Indeed,  let $v\in \cf K_{n}$ and
 $w\in VX\setminus VK_{n}$ be such that $d(v,w)=1$.  Then, there are
 $\g_{0},\g_{1}\in\G$, $v_{0},v_{1}\in VF$, such that $v=\g_{0}
 v_{0}$ and
 $w=\g_{1} v_{1}$.  Moreover, we have $\g_{0}\in E_{n}$ and $\g_{1}\not\in
 E_{n}$.  Let $\d:= \g_{0}^{-1}\g_{1}$, so that $dist(F,\d F) =
 dist(\g_{0}F,\g_{1}F) \leq d(v,w) = 1$, which implies that $\d\in A$.
 Hence, $\g_{0}\in U_{n}$, and the claim follows.

 Finally,
 \begin{align*}
     |\cf K_{n}| &\leq |U_{n}| \cdot |F|\\
     & \leq b \sum_{\d\in A} |E_{n}\setminus E_{n}\cdot \d^{-1}|
     \\
     & = b \sum_{\d\in A} |E_{n}\cdot \d \setminus E_{n}| \\
     & \leq b |A| \cdot |E_{n}\cdot A \setminus E_{n}| \\
     & \leq \frac{b}{n}|E_{n}| \leq \frac{b}{an} |VK_{n}|,
 \end{align*}
 so condition $(iii)$ of Definition \ref{def:amenableGraph} is
 satisfied, showing that $\set{ K_{n}: n\in\bn}$ is an amenable
 exhaustion. Hence, $X$ is amenable, as desired.
 \end{proof}

 If $\Omega \subset VX$, $r\in\bn$, we write $B_{r}(\O) := \set{ v'\in
 VX: \r(v',v)\leq r}$, where $\r$ is the geodesic metric on $VX$.

 \begin{Lemma}
     Let $(X,\G,F)$ be as above.  Let $d:=\sup_{v\in VX} \deg(v)
     <\infty$.  Let $\{K_{n}\}$ be an amenable exhaustion of $X$, and
     $\eps_n := \frac{|\cf K_n|}{|VK_n|} \to 0$.  Then, for any
     $r\in\bn$, $|B_r(\cf K_n)| \leq (d+1)^r \eps_n |VK_n|$.
 \end{Lemma}
 \begin{proof}
     Since
     $$
     B_{r+1}(v)=\bigcup_{v'\in B_{r}(v)}B_{1}(v'),
     $$
     we have $|B_{r+1}(v)|\leq (d+1) |B_{r}(v)|$, giving $|B_{r}(v)|\leq
     (d+1)^{r}$, $\forall v\in VX$, $r\geq 0$.  As a consequence, for any
     finite set $\Omega\subset VX$, we have $B_{r}(\Omega)=\cup_{v'\in
     \Omega}B_{r}(v')$, giving
     \begin{equation}\label{subsets}
     |B_{r}(\Omega)|\leq |\Omega| (d+1)^{r},\quad\forall r\geq 0.
     \end{equation}
     Therefore, $|B_{r}(\cf K_{n})| \leq (d+1)^r |\cf K_n| = (d+1)^r \eps_n
     |VK_n|$.
 \end{proof}

 \begin{Lemma}\label{lem:limitTrace}
     Let $(X,\G,F)$ be as above.  Let $\{K_{n}\}$ be an amenable
     exhaustion of $X$.  Then, for any $B\in\cn_0(X,\G)$, we have
     $$
     \lim_{n\to\infty} \frac{Tr(P(K_n) BP(K_n))}{|VK_{n}|} =
     \frac{1}{|\cf_0|}\ Tr_\G(B),
     $$
     where $P(K_{n})$ is the orthogonal projection of $\ell^{2}(VX)$
     onto $\ell^{2}(VK_{n})$.
 \end{Lemma}
 \begin{proof}
     Denote by $\cf_0$ a subset of $VF$ consisting of one
     representative vertex for each $\G$-class, and let $\cf' :=
     VF\setminus \cf_0$ and $\d:= diam F$.  Then, for any $n\in\bn$,
     $VK_n= \sqcup_{\g\in E_n} \g \cf_0 \sqcup  \O_n$, where
     $\sqcup$ denotes ``disjoint union'' and $\O_n\subset B_{\d}(\cf
     K_n)$.  Indeed, if $v\in \O_n := VK_n \setminus \sqcup_{\g\in
     E_n} \g \cf_0$, then there is a unique $\g\in\G$ such that
     $v\in\g\cf_0$, so that $\g\not\in E_n$, which implies $\g F \cap
     (VX\setminus VK_n) \neq \emptyset$, and $d(v,VX\setminus VK_n)
     \leq \d$, which is the claim.  Therefore,
     \begin{align*}
     Tr(P(K_n)B) & = \sum_{v\in VK_n} B(v,v) \\
     &= \sum_{\g\in E_n} \sum_{v\in\cf_0} B(\g v,\g v) +
     \sum_{v\in\O_n} B(v,v) \\
     &= \sum_{\g\in E_n} \sum_{v\in\cf_0} B(v,v) + \sum_{v\in\O_n}
     B(v,v) \\
     &= |E_n| Tr_\G(B) + \sum_{v\in\O_n} B(v,v).
     \end{align*}
     Moreover,
     $$
     \left| \sum_{v\in\O_n} B(v,v) \right| \leq \|B\| |\O_n| \leq
     \|B\| |B_\d(\cf K_n)| \leq (d+1)^\d \|B\| \eps_n |VK_n|,
     $$
     so that
     $$
     \lim_{n\to\infty} \frac{\sum_{v\in\O_n} B(v,v)}{|VK_n|} = 0.
     $$
     Besides,
     $$
     \lim_{n\to\infty} \frac{|E_n|}{|VK_n|} =
     \frac{1}{|\cf_0|},
     $$
     because $|VK_n| = |E_n| \cdot |\cf_0| + |\O_n|$.  The claim
     follows.
 \end{proof}

 \begin{Lemma}\label{lem:normEstimate}
     Let $(X,\G)$ be as above.  Let $A$ and $Q$ be as in Section
     \ref{sec:DetFormula}.  Let $f(u) := Au-Qu^2$, for $u\in\bc$.
     Then $\|f(u)\|<\frac12$, for $|u|<\frac{1}{d+\sqrt{d^2+2(d-1)}}$.
 \end{Lemma}
 \begin{proof}
     This follows from the estimate
     $$
     \|f(u)\|  \leq |u| \|A\| + |u|^2\|Q\| \leq d|u|+(d-1)|u|^2,
     $$
     which is valid for any $u\in\bc$.
\end{proof}

 \begin{Thm}[Approximation by finite graphs]
     Let $X$ be a connected, countably infinite graph, and let $\G$ be
     a countable discrete amenable subgroup of automorphisms of $X$,
     which acts on $X$ freely and cofinitely, and let $F$ be a
     corresponding fundamental domain.  Let $\set{K_{n} :
     n\in\bn}$ be an amenable exhaustion of $X$.  Then
     $$
     Z_{X,\G}(u) = \lim_{n\to\infty}
     Z_{K_{n}}(u)^{\frac{|\cf_0|}{|K_{n}|}},
     $$
     uniformly on compact subsets of $\set{ u\in\bc: |u|
     <\frac{1}{ d+\sqrt{d^2+2(d-1)} } }$.
 \end{Thm}
 \begin{proof}
     For a finite subset $N\subset VX$, denote by $P(N)\in \cb(\ell^2
     (VX))$ the orthogonal projection of $\ell^{2}(VX)$ onto
     $\text{span} (N)$.  Observe that, since $N$ is an
     orthonormal basis for $\ell^2 (N)$, we have $Tr \bigl(P(N) \bigr) =
     |N|$.

     Let $f(u):= Au-Qu^{2}$ and $P_{n}:= P(VK_{n})$. Then
     $$
     \log Z_{K_{n}}(u) = -\frac12 Tr(P_{n}(Q-I)P_{n})
     \log(1-u^{2}) - Tr \log(P_{n}(I-f(u))P_{n}).
     $$
     Moreover,
     $$
     Tr \log(P_{n}(I-f(u))P_{n}) = -\sum_{k=1}^{\infty}
     \frac{1}{k} Tr\bigl( (P_{n}f(u)P_{n})^{k}\bigr).
     $$
     Observe that, for $k\geq 2$,
     \begin{align*}
     Tr \bigl( P_{n}f(u)^{k}P_{n} \bigr) &= Tr\bigl(
     P_{n} ( f(u)(P_{n}+P_{n}^{\perp}) )^{k} P_{n} \bigr)\\
     &= Tr\bigl( (P_{n}f(u)P_{n})^{k} \bigr) +
     \sum_{ \substack{ \s\in\{-1,1\}^{k-1}\\ \s\neq \{1,1,\ldots,1\} } }
     Tr \bigl( P_{n} \prod_{j=1}^{k-1} [f(u)P_{n}^{\s_{j}}]f(u)P_{n}
     \bigr),
     \end{align*}
     where $P_n^{-1}$ stands for $P_n^\perp$, the projection onto the
     orthogonal complement of $\ell^{2}(VK_{n})$ in $\ell^{2}(VX)$, and
     \begin{align*}
     | Tr \bigl( P_{n} \prod_{j=1}^{k-1} [f(u)P_{n}^{\s_{j}}] f(u)P_{n}
     \bigr) |  & = | Tr \bigl( ...P_{n}  f(u)P_{n} ^\perp... \bigr) | \\
     & \leq \| f(u) \|^{k-1} Tr ( |P_{n}  f(u)P_{n}^\perp|).
     \end{align*}
     Moreover, with $\O_{n}:= B_{1}(VK_{n})\setminus VK_n\subset
     B_{1}(\cf K_{n})$, we have
     \begin{align*}
     Tr ( |P_{n} f(u)P_{n}^\perp|) & = Tr (|P(K_{n} ) f(u)P(\O_n)|
     )\\
     & \leq \|f(u)\| Tr(P(\O_n)) \\
     & = \|f(u)\| |\O_n| \\
     & \leq \|f(u)\| (d+1) \eps_n |VK_n|.
     \end{align*}
     Therefore, we obtain
     $$
     | Tr \bigl( P_{n}f(u)^{k}P_{n} \bigr) - Tr\bigl(
     (P_{n}f(u)P_{n})^{k} \bigr) | \leq (2^{k-1}-1) \| f(u) \|^{k} (d+1)
     \eps_n |VK_n|,
     $$
     so that
     \begin{align*}
     \biggl| Tr&\log(P_{n}(I-f(u))P_{n}) -
     Tr(P_{n} \log(I-f(u)) P_{n}) \biggr| \\
     &= \biggl| \sum_{k=1}^{\infty} \frac{1}{k} Tr\bigl(
     (P_{n}f(u)P_{n})^{k} \bigr) - \sum_{k=1}^{\infty}
     \frac{1}{k} Tr \bigl( P_{n}f(u)^{k}P_{n} \bigr)\biggr| \\
     & \leq \biggl( \sum_{k=1}^{\infty}
     \frac{2^{k-1}\|f(u)\|^{k}}{k}\biggr)(d+1) \eps_n |VK_n| \\
     & \leq C (d+1) \eps_n |VK_n|,
     \end{align*}
     where the series converges for $|u|<
     \frac{1}{ d +\sqrt{d^2+2(d-1)} }$, by Lemma \ref{lem:normEstimate}.
     Hence,
     $$
     \biggl| \frac{Tr\log(P_{n}(I-f(u))P_{n})}{|VK_{n}|} -
     \frac{Tr(P_{n} \log(I-f(u)) P_{n})}{|VK_{n}|} \biggr| \to 0, \
     n\to\infty
     $$
     and, by using Lemma \ref{lem:limitTrace},
     \begin{align*}
     \lim_{n\to\infty} \frac{\log Z_{K_{n}}(u)}{|VK_{n}|} &= -\frac12
     \lim_{n\to\infty} \frac{Tr(P_{n}(Q-I)P_{n}) }{|VK_{n}|}
     \log(1-u^{2}) \\
     & \qquad - \lim_{n\to\infty} \frac{Tr(P_{n} \log(I-f(u))
     P_{n})}{|VK_{n}|} \\
     &= -\frac{1}{|\cf_0|}\left( \frac12 Tr_{\G}(Q-I) \log(1-u^{2}) +
     Tr_{\G}(\log(I-f(u))) \right) \\
     & = \frac{1}{|\cf_0|} \ \log Z_{X,\G}(u),
     \end{align*}
     from which the claim follows.
 \end{proof}

  \begin{rem}
     Observe that $\frac{1}{2\a}<\frac{1}{d+\sqrt{d^2+2(d-1)}}<\frac{1}{\a}$.
 \end{rem}

 \begin{ack}
     The second and third named authors would like to thank
     respectively the University of California, Riverside, and the
     University of Roma ``Tor Vergata'' for their hospitality at
     various stages of the preparation of this paper.
 \end{ack}


\end{document}